# A Discussion on Nonlinear Quadratic Control and Sontag's Formula

Boris Lohmann and Joscha Bongard, Chair of Automatic Control, TUM School of Engineering and Design, Technical University of Munich, Lohmann@tum.de

The quadratic optimal state feedback (LQR) is one of the most popular designs for linear systems and succeeds via the solution of the algebraic Riccati equation. The situation is different in the case of non-linear systems: the Riccati equation is then replaced by the *Hamilton Jacobi Bellman equation* (HJB), the solution of which is generally difficult. A compromise can be the so-called *Inverse Optimal Control*, a form of which is *Sontag's formula* [1]; here the minimized cost function follows from the feedback law chosen, not vice versa. Using Sontag's formula in the variant according to Freeman and Primbs [2, 9], the actually minimized cost function is given in the following sections, including cases when it reduces to the quadratic cost. Also some remarks and thoughts are presented for discussion.

## Sontag's Formula and Optimality

We consider an input affine nonlinear system of order *n* with a vector $u$ of input variables,

$$\dot{x} = f(x) + G(x)u \, , \text{ with } f(0) = 0, \; f, G \in C^1, \; x(0) = x_0. \tag{1}$$

Yet without knowing the feedback law, we can try to say something about the (future) stability behavior of the equilibrium $0$ by considering a positive definite function $V(x)$ as a possible Lyapunov function and by writing down $\dot{V}(x)$:

$$\dot{V} = \left(\frac{\partial V}{\partial x}\right)^T \dot{x} = \underbrace{\left(\frac{\partial V}{\partial x}\right)^T f(x)}_{a(x)} + \underbrace{\left(\frac{\partial V}{\partial x}\right)^T G(x)}_{b^T(x)} u \, . \tag{2}$$

If $b(x) \neq 0$, then by a suitable choice of $u$ we can assure $\dot{V} < 0$, and thus asymptotic stability of the equilibrium $0$. If, on the other hand, $b(x) = 0$, then $a(x) < 0$ must apply for asymptotic stability. This motivates the following definition:

Def.: A differentiable, positive definite and radially unbounded function $V(x)$ is called a *Control Lyapunov Function* (CLF) of the system (1), if for all vectors $x \neq 0$, with which $b(x) = 0$, we will have: $a(x) < 0$.

***Sontag's Formula*** ([1], in the extended version according to [2, 9], see also [3, 5, 8]):

The state feedback

$$u = \begin{cases} -R^{-1}b \cdot \underbrace{\dfrac{a + \sqrt{a^2 + x^T Q x \, b^T R^{-1} b}}{b^T R^{-1} b}}_{\lambda(x)} & \text{for } b(x) \neq 0 \\ 0 & \text{for } b(x) = 0 \end{cases} \qquad (3)$$

(where $Q$, $R$ are chosen constant, symmetric und positive definite and $V(x)$ is a CLF) makes $0$ a globally asymptotically stable equilibrium, at the same time minimizing the cost function

$$J = \frac{1}{2} \int_0^\infty \frac{1}{\lambda(x)} \left( x^T Q x + u^T R u \right) dt \qquad (4)$$

for any trajectory of the closed-loop system.

The cost function (4) is "distorted" by the factor $1/\lambda(x)$ compared to the *classical quadratic cost function*

$$J = \frac{1}{2} \int_0^\infty \left( x^T Q x + u^T R u \right) dt \; . \qquad (5)$$

According to (3), $\lambda(x)$ depends on the selected CLF $V(x)$ as well as on $f(x)$, $G(x)$, $Q$ and $R$. However, since $\lambda > 0$ (for $x \neq 0$) is valid, the cost function (4) remains transparent and interpretable. In particular, a compromise between a "low" actuation effort and a "nice" trajectory can be set by choosing $Q$, $R$ appropriately. The feedback law (3) is continuous in the neighborhood of $b(x) = 0$, as one can easily show with the rule of l'Hospital. A proof of the result (4) can be found in appendix 1 (this proof was given by Michael Buhl in 2006).

## Discussion

- Unfortunately, in most publications a cost function of the type (5) is given in connection with the control law (3), which can lead to misunderstandings. For although the above control design with a *freely specifiable* CLF is *oriented* towards the objective (5), this objective is by far not achieved with *every* CLF, but only with the (hard to find) *solution* $V(x) = J^*(x)$ of the *Hamilton Jacobi Bellman equation* [2, 3, 5, 7]

$$\frac{1}{2} x^T Q x - \frac{1}{2} J_x^{*T} G R^{-1} G^T J_x^* + J_x^{*T} f = 0 \qquad (6)$$

  for the problem (1), (5) with $J^*(0) = 0$. This solution $J^*(x)$ is used in the feedback

$$u = -R^{-1} G^T J_x^* = -R^{-1} b(x), \qquad (7)$$

  where $J^*(x)$ is the *optimal value function* (cost-to-go function) and $J_x^*$ the corresponding derivative for $x$. It turns out that $\lambda(x) \equiv 1$ in (3) (see appendix 2), meaning that Sontag's formula (3) minimizes (5) if the solution $J^*(x)$ of the HJB (6) is used as the CLF.

For *freely* specifiable CLF, it is only in [3] that the actually minimized cost function (4) is clearly given in equation (3.5.28), but with an error that crept in through the line before it[1]. In [5], too, the result (4) can be read on page 94.

- It should be emphasized that in the publications [3, 5, 9, 4, 10, 6, 7] and others one can find numerous more in-depth contributions to *Inverse Optimal Control*, to generalizations of Sontag's formula, its robustness, and to references to the Hamilton Jacobi Bellman equation. Suggestions for the determination of suitable CLFs are also made.
- If we have the choice between different CLFs, we will probably prefer one that makes $\lambda$ "as constant as possible", at least in the vicinity of the equilibrium $\boldsymbol{0}$, because $\lambda = const$ passes (4) into the classical cost function (5).
- For systems of order $n=1$ and by using $V = \frac{1}{2}x^2$, the feedback (3) simultaneously provides the optimal solution in the sense of (5), see [5], page 95 and [3], page 117.

## Appendix 1: Proof of optimality of the feedback (3)

To derive Sontag's formula (3), we start from the problem (1), (4), with still open function $\lambda(\boldsymbol{x}) > 0$ (for $\boldsymbol{x} \neq \boldsymbol{0}$). We use the Hamilton Jacobi Bellman equation in its general form

$$\min_{u} \left\{ \frac{1}{2\lambda(\boldsymbol{x})} \left( \boldsymbol{x}^T \boldsymbol{Q} \boldsymbol{x} + \boldsymbol{u}^T \boldsymbol{R} \boldsymbol{u} \right) + J_{\boldsymbol{x}}^{*T}(\boldsymbol{x}) \left( \boldsymbol{f}(\boldsymbol{x}) + \boldsymbol{G}(\boldsymbol{x}) \boldsymbol{u} \right) \right\} = 0. \tag{8}$$

Since there are no constraints of variables to consider, the minimum search can be done by derivation of the curly bracket with respect to $\boldsymbol{u}$ and setting it to zero:

$$\frac{\partial}{\partial \boldsymbol{u}} \left\{ \frac{1}{2\lambda(\boldsymbol{x})} \left( \boldsymbol{x}^T \boldsymbol{Q} \boldsymbol{x} + \boldsymbol{u}^T \boldsymbol{R} \boldsymbol{u} \right) + J_{\boldsymbol{x}}^{*T}(\boldsymbol{x}) \left( \boldsymbol{f}(\boldsymbol{x}) + \boldsymbol{G}(\boldsymbol{x}) \boldsymbol{u} \right) \right\} = \frac{1}{\lambda} \boldsymbol{R} \boldsymbol{u} + \boldsymbol{G}^T J_{\boldsymbol{x}}^* = \boldsymbol{0} \tag{9}$$

$$\Rightarrow \quad \boldsymbol{u} = -\lambda \boldsymbol{R}^{-1} \boldsymbol{G}^T J_{\boldsymbol{x}}^*. \tag{10}$$

Now we insert the result (10) into (8), obtaining

$$\frac{1}{2\lambda} \left( \boldsymbol{x}^T \boldsymbol{Q} \boldsymbol{x} + \lambda^2 J_{\boldsymbol{x}}^{*T} \boldsymbol{G} \boldsymbol{R}^{-1} \boldsymbol{G}^T J_{\boldsymbol{x}}^* \right) + J_{\boldsymbol{x}}^{*T} \left( \boldsymbol{f} - \boldsymbol{G} \lambda \boldsymbol{R}^{-1} \boldsymbol{G}^T J_{\boldsymbol{x}}^* \right) = 0. \tag{11}$$

Due to $\lambda > 0$, we may multiply (11) by $\lambda$, obtaining the *Hamilton Jacobi Bellman equation* for the problem at hand in the form

$$\frac{1}{2} \boldsymbol{x}^T \boldsymbol{Q} \boldsymbol{x} - \frac{1}{2} \lambda^2 \underbrace{J_{\boldsymbol{x}}^{*T} \boldsymbol{G}}_{\boldsymbol{b}^T} \boldsymbol{R}^{-1} \underbrace{\boldsymbol{G}^T J_{\boldsymbol{x}}^*}_{\boldsymbol{b}} + \lambda \underbrace{J_{\boldsymbol{x}}^{*T} \boldsymbol{f}}_{a} = 0. \tag{12}$$

---

[1] $\dot{V}_{\frac{u_s}{2}}$ is not $-\frac{1}{2} p \boldsymbol{b}^T \boldsymbol{b}$ but is (with $c_0 = 0$) $a - \frac{1}{2} p \boldsymbol{b}^T \boldsymbol{b} = \frac{1}{2}\left(a - \sqrt{a^2 + (\boldsymbol{b}^T \boldsymbol{b})^2}\right) \cdot \frac{a + \sqrt{a^2 + (\boldsymbol{b}^T \boldsymbol{b})^2}}{a + \sqrt{a^2 + (\boldsymbol{b}^T \boldsymbol{b})^2}} = \ldots = -\frac{\boldsymbol{b}^T \boldsymbol{b}}{2p} = -l(\boldsymbol{x})$.

Given the function $\lambda(x)$, we would have to solve the partial differential equation (12) for $J^*$ next. But now we take a *different* route: Even without knowing $J^*$, we can establish a connection between $J^*$ and $\lambda(x)$ by solving the equation (12) for $\lambda(x)$:

$$\begin{cases} \lambda_{1/2} = \dfrac{a \pm \sqrt{a^2 + x^T Q x \, b^T R^{-1} b}}{b^T R^{-1} b} & \text{for } b \neq 0 \\ \lambda = -\dfrac{x^T Q x}{2a} & \text{for } b = 0 \end{cases} \quad (13)$$

For $b \neq 0$, we observe: Since $\lambda$ must be greater than zero, only the first of the two solutions, i.e. $\lambda = \lambda_1$, is suitable. For $b = 0$, $a < 0$ is required[2] in order to ensure $\lambda > 0$. Since $J^*$ is also positive definite, $J^*$ must have all the properties of a *Control Lyapunov Function*:

- $J^*(x)$ is positive definite and
- if $b = G^T J_x^* = 0$, then $a = J_x^{*T} f < 0$ must hold.

Therefore, instead of *calculating* $J^*$ (depending on $\lambda$ and the other parameters of the cost function), we can also *specify* an *arbitrary* Control Lyapunov Function $V(x)$ for the system (1) and immediately give the control law (10) with $J^* = V$:

$$u = -\lambda R^{-1} G^T J_x^* = -R^{-1} b \lambda \qquad \text{as in equation (3).}$$

With $\lambda$ according to (13), the cost function adapts to the choice of $V(x)$. This is also referred to as *Inverse Optimal Control*, because the control law is fixed by the specification of $V(x)$, and the cost function follows that choice.

## Appendix 2: Nonlinear Quadratic Control and Sontag's Formula

Nonlinear quadratic control means minimizing the quadratic cost (5) by the feedback (7) using the solution $J^*(x)$ of the HJB (6) (with $J^*(0) = 0$). Sontag's Formula (3) delivers the *same* feedback, if we use the solution $V(x) = J^*(x)$ of (6) as the CLF in (3). To prove this, we will show that $\lambda(x) \equiv 1$, rendering (3) into (7):

$a(x)$ and $b(x)$, defined in (2), fulfil equation (6),

$$a = \frac{1}{2} b^T R^{-1} b - \frac{1}{2} x^T Q x \;\; \Rightarrow \;\; a^2 = \frac{1}{4}(b^T R^{-1} b)^2 - \frac{1}{2} x^T Q x \, b^T R^{-1} b + \frac{1}{4}(x^T Q x)^2, \quad (14)$$

when using $V(x) = J^*(x)$ as the CLF in Sontag's formula. With (14), we find from (13) for $b \neq 0$:

---

[2] With $a < 0$ and by using l'Hospital's rule it is easy to show that the lower part of (13) is just the continuous continuation of the upper part of (13) at $b = 0$. The idea of solving (12) for $\lambda$ also appears in [9].

$$\lambda(x) = \frac{a + \sqrt{a^2 + x^T Q x \, b^T R^{-1} b}}{b^T R^{-1} b} =$$

$$= \frac{\frac{1}{2} b^T R^{-1} b - \frac{1}{2} x^T Q x + \sqrt{\frac{1}{4}(b^T R^{-1} b)^2 - \frac{1}{2} x^T Q x \, b^T R^{-1} b + \frac{1}{4}(x^T Q x)^2 + x^T Q x \, b^T R^{-1} b}}{b^T R^{-1} b} =$$

$$= \frac{\frac{1}{2} b^T R^{-1} b - \frac{1}{2} x^T Q x + \left| \frac{1}{2} b^T R^{-1} b + \frac{1}{2} x^T Q x \right|}{b^T R^{-1} b} = 1,$$

and by similar steps: $\lambda = 1$ for $b = 0$.

## References

[1]   E. D. Sontag: A 'universal' construction of Artstein's theorem on nonlinear stabilization. Syst. Contr. Lett., vol. 13, no. 2, pp. 117-123, 1989

[2]   R.A. Freeman and J.A. Primbs: Control Lyapunov Functions: New Ideas from an old Source. Proc. of CDC, Kobe, Japan, 1996

[3]   Sepulchre, R., Jankovic, M., Kokotovic, P.V.: Constructive Nonlinear Control. Springer 1997

[4]   J.W. Curtis III: A Generalization of Sontag's Formula for High-Performance CLF-Based Control. Dissertation, Brigham Young University, 2002

[5]   Sackmann, M.: Modifizierte optimale Regelung: Nichtlinearer Reglerentwurf unter Verwendung der Hyperstabilitätstheorie. Dissertation, Univ. Karlsruhe, VDI Fortschritt-Bericht Nr. 906, 2001

[6]   Sackmann, M., Krebs, V.: Modified Optimal Control: Global Asymptotic Stabilization of Nonlinear Systems. Proceedings of IFAC Conf. on Control Systems Design, Bratislava, 2000, p. 199-204.

[7]   J.A. Primbs, V. Nevistic, and J. C. Doyle: Nonlinear Optimal Control: A Control Lyapunov Function and Receding Horizon Perspective. Asian Journal of Control, Vol. 1, No. 1, pp. 14-24, March 1999

[8]   J.A. Primbs and M. Gianelli: Control Lyapunov Function based Receding Horizon Control for Time-Varying Systems. Proceedings of CDC Tampa, Florida, 1998

[9]   J.A. Primbs: Nonlinear Optimal Control: A Receding Horizon Approach. Dissertation, California Institute of Technology, 1999

[10]  R.A. Freeman and P.V. Kokotovic: Robust Nonlinear Control Design. Birkhäuser 1996